\documentclass[11pt]{article}
\usepackage{latexsym}
\title{Curvature Estimates and the Positive Mass Theorem}
\author{Hubert Bray\thanks{Research supported by the NSF, Grant No.\
DMS-9706006.}\,
and Felix Finster}
\date{June 1999}

\newtheorem{Def}{Def.}[section]
\newtheorem{Thm}[Def]{Theorem}

\newtheorem{Lemma}[Def]{Lemma}
\newtheorem{Corollary}[Def]{Corollary}

\newcommand{\Proof}{{\em{Proof: }}}
\newcommand{\QED}{\ \hfill $\FBox$ \\[1em]}

\newcommand{\spc}{\;\;\;\;\;\;\;\;\;\;}
\newcommand{\bra}{\mbox{$< \!\!$ \nolinebreak}}
\newcommand{\ket}{\mbox{\nolinebreak $>$}}

\newcommand{\R}{\mbox{\rm I \hspace{-.8 em} R}}
\newcommand{\1}{\mbox{\rm 1 \hspace{-1.05 em} 1}}

\newcommand{\Sl}{\mbox{$\prec \!\!$ \nolinebreak}}
\newcommand{\Sr}{\mbox{\nolinebreak $\succ$}}

\newcommand{\FBox}{\rule{2mm}{2.25mm}}

\begin{document}
\maketitle

\begin{abstract}
The Positive Mass Theorem implies that any smooth, complete, 
asymptotically flat 3-manifold with non-negative scalar curvature which has
zero total mass is isometric to $(\R^3, \delta_{ij})$.  In this paper,
we quantify this statement using spinors and prove that if a complete,
asymptotically flat manifold with non-negative scalar curvature has small
mass and bounded isoperimetric constant, then the 
manifold must be close to $(\R^3,\delta_{ij})$, in the sense that there is 
an upper bound for the $L^2$ norm of the Riemannian curvature tensor over
the manifold except for a set of small measure. 
This curvature estimate allows us 
to extend the case of equality of the Positive Mass Theorem to include
non-smooth manifolds with generalized non-negative scalar curvature, which we
define.  
\end{abstract}

\section{Introduction}
We introduce our problem in the context of General Relativity. Consider a
$3+1$ dimensional Lorentzian manifold $N$ with metric $g_{\alpha \beta}$ of
signature $(-+++)$.
We denote the induced Levi-Civita connection by $\bar{\nabla}$. Then
the corresponding
Ricci tensor $\bar{R}_{\alpha \beta}$ satisfies Einstein's equations
\begin{equation}
\bar{R}_{\alpha \beta} - \frac{1}{2} \:\bar{R}\:g_{\alpha \beta} \;=\; 8 \pi\:
T_{\alpha
\beta} \;\;\;, \label{00}
\end{equation}
where $T_{\alpha \beta}$ is the energy-momentum tensor (which describes
the distribution of matter in space-time). Furthermore, we are given a
complete,
oriented, space-like hypersurface $M$. The Lorentzian metric $g_{\alpha
\beta}$ induces on $M$ a Riemannian metric
$g_{ij}$ (we always use Latin indices on the hypersurface and Greek
indices in the embedding manifold). Choosing on $M$ a normal vector
field $\nu$, the exterior curvature of $M$ is given by
the second fundamental form $h_{jk} = (\bar{\nabla}_j \nu)_k$.
We make the following assumptions:
\begin{description}
\item[(I)] $M$ is asymptotically flat. Thus
we assume that there is a compact set $K$ such
that $M \setminus K$ is diffeomorphic to the region $\R^3 \setminus B_r(0)$
outside a ball of radius $r$. Under this diffeomorphism, the
metric should be of the form
\begin{equation}
g_{jk}(x) \;=\; \delta_{jk} + a_{jk}(x) \;\;\;,\spc x \in \R^3 \setminus
B_r(0) \;\;\;,
\label{1}
\end{equation}
where $a_{ij}$ decays at infinity as
\begin{equation}
a_{ij} \;=\; O(1/r) \;\;\;,\;\;\;\;\;
\partial_k a_{ij} \;=\; O(1/r^2) \;\;\;{\mbox{, and}} \;\;\;\;\;
\partial_{kl} a_{ij} = O(1/r^3) \;\;\;. \label{2c}
\end{equation}
The second fundamental form should decay as
\begin{equation}
h_{ij} = O(1/r^2) \spc {\mbox{and}} \spc \partial_k h_{ij} = O(1/r^3)
\;\;\;.
\label{3a}
\end{equation}
\item[(II)] The energy-momentum tensor satisfies on $M$ the dominant
energy
condition \cite{HE}, i.e.\ for each point $p \in M$ and for each
time-like
vector $u \in T_pN$, the vector $T^\alpha_\beta u^\beta$ is
non-spacelike and $T_{\alpha \beta} \:u^\alpha u^\beta \leq 0$.
\end{description}
For the hypersurface $M$, one can define the total energy and
momentum, as first introduced by Arnowitt, Deser, and Misner \cite{ADM}.
For this, one considers in the asymptotic end the coordinate
spheres $S_R$, $R>r$, around
the origin and takes limits of integrals over these spheres,
\begin{eqnarray}
E &=& \frac{1}{16 \pi} \lim_{R \rightarrow \infty} \int_{S_R}
(\partial_j g_{ij} - \partial_i g_{jj}) \:d\Omega^i \label{12} \\
P_k &=& \frac{1}{8 \pi} \lim_{R \rightarrow \infty} \int_{S_R}
(h_{ki} - \delta_{ki} \:h_{jj}) \:d\Omega^i \;\;\; , \label{13}
\end{eqnarray}
where $d\Omega^i = \nu^i \:du$, $du$ is the area form, and $\nu$ is
the normal vector to $S_R$ in the coordinate chart.
The Positive Energy Theorem \cite{SY1, W} states that, under the
considered assumptions, $E \geq |P|$. In the case that the second fundamental
form is identically zero, the total momentum vanishes. Then the total
energy is also called the total mass $m$, and the Positive Mass Theorem
states that $m \geq 0$.

In this paper, our aim is to study how total energy and momentum
control the Riemannian curvature tensor. Following Witten's proof of the
Positive Energy Theorem, we consider the massless Dirac equation on the
hypersurface $M$. We derive an integral estimate for the Riemannian curvature
tensor $\bar{R}_{\alpha \beta \gamma \delta}$ involving total
energy/momentum and the Dirac wave function $\Psi$.
We then restrict our attention to the case of zero
second fundamental form. By substituting a-priori estimates for the Dirac
wave function, we get an $L^2$ estimate for the Riemannian curvature
tensor of $M$ in terms of its total mass.
More precisely, our main theorem is the following:
\begin{Thm}
\label{thm6.2} There exist positive constants $c_1$, $c_2$, and $c_3$
such that for
any smooth, complete, asymptotically flat manifold $(M^3,g)$ with non-negative
scalar curvature and total mass $m$ and any smooth, bounded function $\eta$
with bounded gradient on $M$, 
\begin{equation}
\int_{M \backslash D} \eta \: R_{ijkl} R^{ijkl} \: d\mu  
\;\leq\;  m \:c_1\:\sup_M \left(|\eta| |R_{ijkl}| + 
|\Delta \eta| \right) 
+ \sqrt{m} \:c_2\: \left\| \eta\: |\nabla_k R_{ij \alpha \beta}|
\right\|_{L^2} \;\;\;,
\end{equation}
where the set $D$ depends on $M$ with
\begin{equation}
   \mbox{Vol}(D)^{1/3} \;\leq\; 64 \pi \:c_3\: \frac{m}{k^2} \;\;\;,
\end{equation}
$k = \inf \frac{A}{V^{2/3}}$ is the isoperimetric constant of $M$,
$R_{ijkl}$ is the Riemannian curvature tensor of $M$, and 
$\|.\|_{L^2}$ is the $L^2$-norm on $(M^3, g)$.
\end{Thm}
As an application of this theorem, we finally extend the case of equality of
the Positive Mass Theorem to non-smooth manifolds.

\section{Spinors, the Hypersurface Dirac Operator}
We begin with a brief introduction to Dirac spinors in curved space-time.
Following~\cite{F}, we work in a coordinate chart (for a coordinate-free
formulation see e.g.\ \cite{PT}). The {\em{Dirac operator}} $G$ is a
differential operator of first order,
\begin{eqnarray}
G \;=\; i G^\alpha(x) \frac{\partial}{\partial x^\alpha} + B(x) \;\;\;,
\label{0}
\end{eqnarray}
where $B$ and the Dirac matrices $G^\alpha$ are $(4 \times 4)$-matrices.
The Dirac matrices and the Lorentzian metric are related by the
anti-commutation relations
\begin{eqnarray}
-g^{\alpha \beta}(x) \;=\; \frac{1}{2} \left\{ G^\alpha(x),\:G^\beta(x)
\right\} \;\equiv\; \frac{1}{2} \left( G^\alpha G^\beta + G^\beta
G^\alpha \right)\!(x) \;\;\; .
\label{4}
\end{eqnarray}
The four-component, complex wave function $\Psi$ of a Dirac particle
satisfies the Dirac equation
\[ (G-m_0) \:\Psi \;=\; 0 \;\;\; , \]
where $m_0$ is the rest mass of the Dirac particle.
At every space-time point $x$, the wave functions are endowed with an
indefinite scalar product, which we call {\em{spin scalar product}}.
For two wave functions $\Psi$ and $\Phi$, it takes the form
\[ \Sl \Psi \:|\: \Phi \Sr(x) \;=\; \Psi^*(x) \left( \begin{array}{cc}
\1 & 0 \\ 0 & -\1 \end{array} \right) \Phi(x) \;\;\;, \]
where `$^*$' denotes complex conjugation and where $\1$, $0$ are $(2
\times 2)$-submatrices (in physics, this scalar product is currently
written in the form $\Sl \Psi \:|\: \Phi \Sr = \overline{\Psi} \Phi$
with the ``adjoint spinor'' $\overline{\Psi}$). The Dirac matrices
$G^\alpha(x)$ are Hermitian with respect to the spin scalar product. By
integrating the spin scalar product over the hypersurface $M$, we form
the scalar product
\begin{eqnarray}
\bra \Psi \:|\: \Phi \ket \;=\; \int_M \Sl \Psi \:|\: G^\alpha \:\Phi
\Sr
\:\nu_\alpha\: d\mu \;\;\; , \label{3}
\end{eqnarray}
where $d\mu=\sqrt{{\mbox{det }} g_{ij}} \:d^3x$ is the invariant measure on
$M$.
This scalar product is definite; we can assume it to be positive.
The integrand of
(\ref{3}) has the physical interpretation as the probability density of
the
particle. Since it will appear in our calculations very often, we
introduce
the short notation
\[ (\Psi \:|\: \Phi) \;\equiv\; \Sl \Psi \:|\: G^\alpha \nu_\alpha \:
\Phi \Sr \;\;\; . \]

For a given Lorentzian metric, the Dirac matrices $G^\alpha(x)$ are
not uniquely determined by the anti-commutation relations (\ref{4}). One
way to handle this problem is to work with spin bundles and
orthonormal frames \cite{PT}. More generally, it is shown in \cite{F}
that
all possible choices of Dirac matrices lead to unitarily equivalent
Dirac
operators.
One must keep in mind, however, that the matrix $B(x)$ in (\ref{0})
depends
on the choice of the $G^\alpha$; it can be given explicitly in terms
of the Dirac matrices $G^\alpha$ and their first partial derivatives.

The Dirac matrices induce a connection $D$ on the spinors, which we
call {\em{spin derivative}}. In a chart, it has the representation
\begin{equation}
        D_\alpha \;=\; \frac{\partial}{\partial x^\alpha} \:-\: i
E_\alpha(x)
        \;\;\; ,
        \label{5}
\end{equation}
where the $(4 \times 4)$-matrices $E_\alpha(x)$ are functions of 
$G^\alpha(x)$ and $\partial_\alpha G^\beta(x)$ (see \cite{F} for an explicit
formula). The spin derivative is compatible with the spin scalar product, i.e.
\begin{equation}
\partial_j \Sl \Psi \:|\: \Phi \Sr \;=\; \Sl D_j \Psi \:|\: \Phi \Sr
\:+\: \Sl \Psi \:|\: D_j \Phi \Sr \;\;\;.
        \label{8}
\end{equation}
Furthermore, the combined spin and covariant derivative of the Dirac
matrices vanishes,
\begin{equation}
[D_j, G^k] \:+\: \bar{\Gamma}^k_{jl} \:G^l \;=\; 0
        \label{7}
\end{equation}
($\bar{\Gamma}^i_{jk}$ denote the Christoffel symbols of the Levi-Civita
connection).
The curvature of the spin connection is
given by the commutator $[D_j,D_k] \equiv D_j D_k - D_k D_j$. It
is related with the Riemannian curvature tensor by
\begin{equation}
[D_\alpha, D_\beta] \;=\; \frac{1}{8} \:
{\bar{R}}_{\alpha \beta \gamma \delta} \: [G^\gamma, G^\delta] \;\;\; .
\label{19a}
\end{equation}
Using the spin derivative, one can write the Dirac operator (\ref{0}) in
the alternative form
\[ G \;=\; i \:\sum_{\alpha=1}^4 G^\alpha(x) \:D_\alpha \;\;\; . \]

Similar to the procedure in \cite{W,PT}, we next define the
so-called {\em{hypersurface Dirac operator}} ${\cal{D}}$. For this, we
consider the Dirac matrices $G^\alpha$ and the spin derivative (\ref{5}) of the
Lorentzian manifold, but take only the derivatives tangential to $M$,
\[ {\cal{D}} \;:=\; i \sum_{j=1}^3 G^j(x) \:D_j \;\;\; . \]
The hypersurface Dirac operator can be considered as a differential
operator
on the four-component wave functions on the hypersurface $M$. According
to \cite{W,PT}, the square of the hypersurface Dirac operator satisfies
the Weitzenb\"ock formula
\begin{eqnarray}
{\cal{D}}^2 &=& D^*_j D^j\:+\: {\cal{R}} \spc {\mbox{with}} 
\label{10} \\
{\cal{R}} &=& \frac{1}{4} \:(\bar{R} \:+\: 2 \bar{R}_{\alpha \beta} \:\nu^\alpha
\nu^\beta \:+\: 2 \bar{R}_{\alpha i} \:\nu^\alpha\:\nu_\beta\:G^\beta G^i)
\;\;\; ,
\nonumber
\end{eqnarray}
where $D_j^*$ denotes the formal adjoint of the operator $D_j$.
As a consequence of the dominant energy condition {\bf{(II)}}
and Einstein's equations (\ref{00}), the $(4 \times 4)$-matrix
${\cal{R}}$ is positive definite.

Let us introduce a convenient notation for the covariant and spin
derivatives. The Levi-Civita connection $\bar{\nabla}_\alpha$ and the spin
connection $D_\alpha$ give a parallel transport of tensor and spinor
fields, respectively. Furthermore, the induced metric $g_{jk}$
yields a Levi-Civita connection on $M$. For clarity, we denote this
last connection and all its derived ``intrinsic'' curvature objects
without a bar; i.e., we have the covariant derivative
${\nabla}_j$ with Christoffel symbols ${\Gamma}^j_{kl}$, the
curvature tensor ${R}^j_{klm}$, etc.. For a derivative tangential 
to $M$, the connections $\bar{\nabla}$ and ${\nabla}$ are related to each 
other by
\begin{equation}
        \bar{\nabla}_j u \;=\; {\nabla}_j u \:+\: h_{jk} \:u^k \:\nu \;\;\; ,
        \label{e}
\end{equation}
where $u$ denotes a vector field tangential to
the hypersurface. Combining the spin and
Levi-Civita connections, we can differentiate all objects with spin
and tensor indices. With a slight abuse of notation, we write this
derivative with a nabla. A bar indicates that we treat the tensor
indices with the Christoffel symbols $\bar{\Gamma}$; otherwise, the
connection ${\Gamma}$ is used. For example, we have
\begin{eqnarray*}
\bar{\nabla}_i \Psi &=& {\nabla}_i \Psi \;=\; D_i \Psi \\
{\nabla}_i {\nabla}_j \Psi &=& D_i D_j \Psi \:-\:
{\Gamma}^k_{ij} \:D_k \Psi \\
\bar{\nabla}_i \bar{\nabla}_j \Psi &=& D_i D_j \Psi \:-\:
\bar{\Gamma}^k_{ij} \:D_k \Psi \;\;\;,\;{\mbox{etc.}}
\end{eqnarray*}

For the proof of our curvature estimates, we shall choose a constant
spinor $\Psi_0$ in the asymptotic end and consider the solution of the
massless hypersurface Dirac equation with asymptotic boundary values
$\Psi_0$,
\begin{equation}
{\cal{D}} \Psi \;=\; 0 \spc{\mbox{with}} \spc \lim_{|x| \rightarrow
\infty} \Psi(x) \;=\; \Psi_0 \;\;\; .
        \label{6}
\end{equation}
The existence of such a solution is proved in \cite{PT}. The wave
function behaves at infinity like
\begin{equation}
\partial_j \Psi \;=\; O(1/r^2) \;\;\;,\spc \partial_{jk} \Psi
\;=\; O(1/r^3) \;\;\; . \label{13c}
\end{equation}
We remark that the massless Dirac equation (\ref{6}) decouples into
two two-spinor equations, the so-called Weyl equations (which
separately describe the left and right handed components of the Dirac
spinor). But this is not very useful for us; we prefer
working with four-component Dirac spinors.

In order to illustrate our notation, we finally outline Witten's proof
of the Positive Energy Theorem in our setting. We take the solution 
$\Psi$ of the Dirac equation (\ref{6}) and compute the following 
divergence:
\begin{eqnarray}
\lefteqn{ {\nabla}_j (D^j \Psi \:|\: \Psi) \;=\;
{\nabla}_j \Sl D^j \Psi \:|\: G^\alpha \nu_\alpha \:\Psi \Sr } 
\nonumber \\
&\stackrel{(\ref{8})}{=}&
\Sl {\nabla}_j D^j \Psi \:|\: G^\alpha \nu_\alpha \:\Psi \Sr \;+\;
\Sl D^j \Psi \:|\: \partial_j(G^\alpha \nu_\alpha) \:\Psi \Sr \;+\;
\Sl D^j \Psi \:|\: G^\alpha \nu_\alpha \:D_j \Psi \Sr \nonumber \\
&\stackrel{(\ref{7})}{=}&
({\nabla}_j D^j \Psi \:|\: \Psi ) \;+\;
\Sl D^j \Psi \:|\: h_{jk} G^k \:\Psi \Sr \;+\;
( D^j \Psi \:|\: D_j \Psi) \nonumber \\
&=& (({\nabla}_j+G^\alpha \nu_\alpha \:h_{jk} G^k) \:D^j \Psi 
\:|\: \Psi) \:+\: (D^j \Psi \:|\: D_j \Psi) \;\;\; .
        \label{a}
\end{eqnarray}
Using that the formal adjoints of the operators $D_j$ are
\begin{equation}
D_j^* \;=\; -{\nabla}_j-G^\alpha \nu_\alpha \:h_{jk} G^k \;\;\;,
        \label{aa}
\end{equation}
we can write (\ref{a}) in the shorter form
\[ (D^j \Psi \:|\: D_j \Psi) \;=\; {\nabla}_j (D^j \Psi \:|\: 
\Psi) \:+\: (D_j^* D^j \Psi \:|\: \Psi) \;\;\; . \]
We now integrate both sides and substitute the Weitzenb\"ock 
formula (\ref{10}). This gives the identity
\begin{equation}
\bra D \Psi \:|\: D \Psi \ket \:+\: \bra \Psi \:|\: {\cal{R}}
\:\Psi \ket \;=\;
\int_M {\nabla}_j (D^j \Psi \:|\: \Psi) \: d\mu \;\;\; .
        \label{13a}
\end{equation}
Since ${\cal{R}}$ is a positive matrix, the left side of (\ref{13a}) is
positive. The right side of this equation is an integral over a
divergence. If this integral is approximated by integrals over the
balls $B_R$, $R>r$, we can apply Gauss' theorem to rewrite them in terms
of boundary integrals over the spheres $S_R$. As explained in detail in
\cite{PT},
these boundary integrals can be identified with the integrals in
(\ref{12}) and (\ref{13}). More precisely,
\begin{eqnarray*}
\int_M {\nabla}_j (D^j \Psi \:|\: \Psi) \:
d\mu &=& \lim_{R \rightarrow \infty} \int_{S_R} (D_j \Psi \:|\: \Psi)\;
d\Omega^j \\
&=& 4 \pi \:(E \:|\Psi_0|^2 + \Sl \Psi_0 \:|\: P_k G^k \:\Psi_0 \Sr )
\;\;\;.
\end{eqnarray*}
The Positive Energy Theorem follows by choosing $\Psi_0$ such that
$\Sl \Psi_0 \:|\: P_k G^k \:\Psi_0 \Sr = -|P|$ and $|\Psi_0|^2=1$.
Namely, in this case, one gets in combination with (\ref{13a}) the
inequalities
\begin{eqnarray}
0 \;\leq\; \bra D\Psi \:|\: D\Psi \ket \;\leq\; 4 \pi \:(E - |P|) \;\;\; .
\label{14b}
\end{eqnarray}

\section{Estimates of the Riemann Tensor}
We begin with a pointwise estimate of the Riemann tensor of the 
Lorentzian manifold in terms of the second derivative of the Dirac wave
function.
\begin{Lemma}
\label{lemma1}
For any solution $\Psi$ of the hypersurface Dirac operator (\ref{6}),
\begin{equation}
\left( \sqrt{\bar{R}_{ijkl} \bar{R}^{ijkl}} - \sqrt{2 \:\bar{R}_{ijk \alpha} 
\nu^\alpha
\:\bar{R}^{ijk \beta} \nu_\beta}
\right)^2 \; (\Psi \:|\: \Psi) \;\leq\; 32 \:
({\nabla}_j {\nabla}_k \Psi \:|\: {\nabla}^j {\nabla}^k
\Psi) \;\;\; . \label{17}
\end{equation}
\end{Lemma}
{\Proof}
Relation (\ref{19a}) and the Schwarz inequality yield the following
estimate:
\begin{eqnarray}
\lefteqn{ (\bar{R}_{jk \alpha \beta} \:G^\alpha G^\beta \Psi \:|\: \bar{R}^{jk
\gamma
\delta} \:G_\gamma G_\delta \Psi) \;=\; 16 \: ([D_j, D_k] \Psi \:|\:
[D^j, D^k] \Psi) } \label{20a} \\
&=& 32 \left( ({\nabla}_j {\nabla}_k \Psi \:|\:
{\nabla}^j {\nabla}^k \Psi) \:-\:
({\nabla}_j {\nabla}_k \Psi \:|\:
{\nabla}^k {\nabla}^j \Psi) \right)
\;\leq\; 64 \:({\nabla}_j {\nabla}_k \Psi \:|\:
{\nabla}^j {\nabla}^k \Psi) \label{20b}
\end{eqnarray}
Let us analyze the curvature term on the left side of (\ref{20a})
more explicitly. For simplicity in notation, we choose a chart with $\nu
=\frac{\partial}{\partial x^0}$. We decompose the Riemann tensor into
the tangential and normal components,
\[ \bar{R}_{jk \alpha \beta} \:G^\alpha G^\beta \;=\; \bar{R}_{jklm} \:G^l G^m
\:+\: 2 \bar{R}_{jkl0} \:G^l G^0 \;\;\; . \]
Since the Dirac matrices are Hermitian with respect to the spin
scalar product and the matrix $G^0$ anti-commutes with the $G^j$,
we obtain
\begin{eqnarray}
\lefteqn{ (\bar{R}_{ij \alpha \beta} \:G^\alpha G^\beta \Psi \:|\: \bar{R}^{ij
\gamma \delta} \:G_\gamma G_\delta \Psi) \;=\;
\Sl \bar{R}_{ij \alpha \beta} \:G^\alpha G^\beta \Psi \:|\: (-G^0)\: \bar{R}^{ij
\gamma \delta} \:G_\gamma G_\delta \Psi \Sr } \nonumber \\
&=& (\Psi \:|\: \bar{R}_{ijkl} \:G^l G^k \:\bar{R}^{ijmn} \:G_m G_n \Psi) \:-\:
4 \:(\Psi \:|\: \bar{R}_{ijk0} \:G^0 G^k \:\bar{R}^{ijm0} \:G_m G_0 \Psi) 
\label{15}
\\
&& -2 \:(\Psi \:|\: \bar{R}_{ijk0} \:G^0 G^k \:\bar{R}^{ijmn} \:G_m G_n \Psi)
 \:+\:2 \:(\Psi \:|\: \bar{R}_{ijkl} \:G^l G^k \:\bar{R}^{ijm0} \:G_m G_0 \Psi)
\;\;\; .\label{16}
\end{eqnarray}
The products of Dirac matrices can be simplified with the
anti-commutation rules (\ref{4}). The important point is that, in both
summands
in (\ref{15}), the Dirac matrices combine to a positive multiple of
the identity,
\begin{eqnarray*}
\bar{R}_{ijkl} \:G^l G^k \:\bar{R}^{ijmn}\: G_m G_n &=& 2 \:\bar{R}_{ijkl} 
\bar{R}^{ijkl}
\:\1 \\
-4 \:\bar{R}_{ijk0} \:G^0 G^k \:\bar{R}^{ijm0} \:G_m G_0 &=& 4 
\:\bar{R}_{ijk\alpha}
\nu^\alpha \: \bar{R}^{ijk \beta} \nu_\beta \:\1 \;\;\; .
\end{eqnarray*}
In the two summands in (\ref{16}), the products of Dirac matrices is
more complicated, and the sign of the terms is undetermined.
But we can bound both summands from below with the Schwarz inequality,
\begin{eqnarray*}
-2 \:(\Psi \:|\: \bar{R}_{ijk0} \:G^0 G^k \:\bar{R}^{ijmn} \:G_m G_n \Psi) &\geq&
- \sqrt{2 \bar{R}_{ijkl} \bar{R}^{ijkl}} \:\sqrt{4 \:\bar{R}_{ijk\alpha}
\nu^\alpha \: \bar{R}^{ijk \beta} \nu_\beta} \: (\Psi \:|\: \Psi) \\
2 \:(\Psi \:|\: \bar{R}_{ijkl} \:G^l G^k \:\bar{R}^{ijm0} \:G_m G_0 \Psi) &\geq&
- \sqrt{2 \bar{R}_{ijkl} \bar{R}^{ijkl}} \:\sqrt{4 \:\bar{R}_{ijk\alpha}
\nu^\alpha \: \bar{R}^{ijk \beta} \nu_\beta} \: (\Psi \:|\: \Psi) \;\;\; .
\end{eqnarray*}
By substituting into (\ref{15}) and (\ref{16}), we obtain the estimate
\[ (\bar{R}_{ij \alpha \beta} \:G^\alpha G^\beta \Psi \:|\: \bar{R}^{ij
\gamma \delta} \:G_\gamma G_\delta \Psi) \;\geq\;
\left( \sqrt{2 \bar{R}_{ijkl} \bar{R}^{ijkl}} \:-\: \sqrt{4 \:\bar{R}_{ijk\alpha}
\nu^\alpha \: \bar{R}^{ijk \beta} \nu_\beta} \right)^2 \:(\Psi \:|\: \Psi)
\;\;\;. \]
\QED
In the following lemma, we estimate the integral over the right side
of (\ref{17}) from above. Similar as in Witten's proof of the Positive Energy
Theorem, this is done by integrating one spin derivative by parts.
The higher order of the derivatives makes the calculation more
complicated; on the other hand, we do not get boundary terms.
\begin{Lemma}
\label{lemma2}
There are constants $c_1$ and $c_2$ independent of the geometry of 
$M$ and $N$ such that for any smooth, bounded function $\eta$
with bounded gradient on $M$ and the 
Dirac wave function of the Positive Mass Theorem (\ref{14b}),
\begin{eqnarray*}
\lefteqn{ \int_M \eta \:({\nabla}^j {\nabla}^k \Psi \:|\: 
{\nabla}_j {\nabla}_k \Psi) } \\
& \leq & c_1 \:(E-|P|) \:\sup_M \left( |\partial_j \eta|
\:|h_{kl}| \:+\: |\eta| \:(|{\nabla}_j h_{kl}| + |R_{ijkl}| +
|h_{ij}|^2) \:+\: |\Delta \eta| \right) \\
&&+ c_2 \:\sqrt{E-|P|} \:\left\|\eta \:(|{\nabla}_k \bar{R}_{ij \alpha
\beta}| \:+\: |h_{ij} \:{\nabla}_k h_{lm}| \:+\: |h_{ij}
\:\bar{R}_{klmn}|) \right\|_{L^2} \; \sup_M |\Psi| \;\;\;.
\end{eqnarray*}
In this formula, $\Delta$ denotes the Laplace-Beltrami operator on $M$,
and $|\Psi| \equiv (\Psi | \Psi)^{\frac{1}{2}}$.
\end{Lemma}
{\Proof}
Exactly as in (\ref{a}), we compute the following divergence:
\[ {\nabla}_j ({\nabla}^j {\nabla}^k \Psi \:|\: {\nabla}_k \Psi) \;=\;
(({\nabla}_j+G^\alpha \nu_\alpha \:h_{jk} G^k)
\:{\nabla}^j {\nabla}^k \Psi \:|\: {\nabla}_k \Psi)
\:+\: ({\nabla}^j {\nabla}^k \Psi \:|\:
{\nabla}_j {\nabla}_k \Psi) \]
Using the short notation with the formal adjoint
\begin{equation}
{\nabla}_j^* \;\equiv\; -{\nabla}_j \:-\: G^\alpha \nu_\alpha 
\: h_{jk} G^k \;\;\;,
        \label{d}
\end{equation}
we can also write
\[ ({\nabla}^j {\nabla}^k \Psi \:|\: {\nabla}_j {\nabla}_k \Psi)
\;=\; {\nabla}_j ({\nabla}^j {\nabla}^k \Psi \:|\:
{\nabla}_k \Psi) \;+\; ({\nabla}_j^* \:{\nabla}^j {\nabla}^k
\Psi \:|\: {\nabla}_k \Psi) \;\;\; . \]
We multiply this equation by $\eta$ and integrate over $M$. 
According to the decay properties (\ref{2c}), (\ref{3a}), and (\ref{13c})
we can integrate by parts without boundary terms and obtain
\begin{eqnarray}
\lefteqn{ \int_M \eta \:({\nabla}^j {\nabla}^k \Psi \:|\:
{\nabla}_j {\nabla}_k \Psi) \:d\mu } \nonumber \\
&=&
-\int_M (\partial_j \eta) \:({\nabla}^j {\nabla}^k \Psi
\:|\: {\nabla}_k \Psi) \:d\mu \;+\;
\int_M \eta \:({\nabla}_j^* {\nabla}^j {\nabla}^k \Psi \:|\:
{\nabla}_k \Psi) \:d\mu \;\;\; . \label{b}
\end{eqnarray}

We estimate the resulting integrals. Since the
left side of (\ref{b}) is real, we must only consider the real
parts of all terms. In the first integral on the right side of (\ref{b}), we can
integrate by parts once again. Using again the decay properties
(\ref{2c}), (\ref{3a}), and (\ref{13c}), and the fact that $\partial \eta$
is bounded, we obtain, as in (\ref{a}),
\begin{eqnarray*}
\lefteqn{ {\mbox{Re }} \int_M (\partial_j \eta) \:({\nabla}^j
{\nabla}^k \Psi \:|\: {\nabla}_k \Psi) \:d\mu } \nonumber \\
&=&\frac{1}{2} \int_M (\partial_j \eta) \left( \partial^j (D^k \Psi
\:|\: D_k \Psi) \:-\: \Sl D^k \Psi \:|\: h^{jl} G_l D_k \Psi \Sr
\right) d\mu \\
&=& \frac{1}{2} \int_M (D^k \Psi \:| \left(-\Delta \eta \:-\:
G^\alpha \nu_\alpha \:G_l \:h^{lj} (\partial_j \eta) \right) D_k \Psi)
\:d\mu \;\;\; .
\end{eqnarray*}
We bound the obtained integral with the sup-norm and substitute the Positive 
Energy Theorem (\ref{14b}),
\begin{eqnarray*}
\left| {\mbox{Re }} \int_M (\partial_j \eta) \:({\nabla}^j
{\nabla}^k \Psi \:|\: {\nabla}_k \Psi) \:d\mu \right|
&\leq& \frac{1}{2} \:\sup_M \left( |\Delta \eta| \:+\: |\partial_j \eta| 
\:|h_{kl}|\right) \:\int_M |D \Psi|^2 \:d\mu \\
&\leq& 2 \pi \:\sup_M \left( |\Delta \eta| \:+\: |\partial_j \eta| 
\:|h_{kl}| \right) \:(E-|P|) \;\;\; .
\end{eqnarray*}

The second summand on the right side of (\ref{b}) is more difficult
because it involves third derivatives of the wave function. Our 
method is to commute the ${\nabla}^k$-derivative to the very left using the 
transformation
\begin{equation}
{\nabla}_j^* {\nabla}^j {\nabla}^k \Psi \;=\;
{\nabla}_j^* \left[ {\nabla}^j, {\nabla}^k \right] \Psi
\:+\: \left[{\nabla}_j^*, {\nabla}^k \right] {\nabla}^j \Psi
\:+\: {\nabla}^k \:{\nabla}_j^* {\nabla}^j \Psi \;\;\; .
        \label{c}
\end{equation}
In the resulting third order term, we can apply the Weitzenb\"ock formula,
\begin{equation}
{\nabla}^k \:{\nabla}_j^* {\nabla}^j \Psi \;=\;
{\nabla}^k (D_j^* D^j) \Psi \;=\; -{\nabla}^k({\cal{R}} 
\:\Psi) \;=\; -(D^k {\cal{R}}) \:\Psi \:-\: {\cal{R}} \:(D^k \Psi) 
\;\;\; . \label{hn}
\end{equation}
The two commutators in (\ref{c}) yield terms involving curvature and 
the second fundamental form, more precisely
\begin{eqnarray}
{\nabla}_j^* \left[ {\nabla}^j, {\nabla}^k \right] \Psi & = & 
\frac{1}{4} \:{\nabla}_j^* \left( \bar{R}^{jk \alpha \beta} \:G_\alpha 
G_\beta \:\Psi \right) \label{f}  \\
\left[{\nabla}_j^*, {\nabla}^k \right] {\nabla}^j \Psi
&\stackrel{(\ref{d})}{=}&
-\left[{\nabla}_j, {\nabla}^k \right] {\nabla}^j \Psi
\:-\: \left[G^\alpha \nu_\alpha h_{jl} G^l,
{\nabla}^k \right] {\nabla}^j \Psi \\
&=&-\frac{1}{4} \:\bar{R}^{jk \alpha \beta} \:G_\alpha G_\beta \:D_j \Psi 
\:-\: {R}^{kl} \:D_l \Psi \:+\: G^\alpha \nu_\alpha 
\:({\nabla}^k h_{jl}) \:G^l D^j \Psi \;\;\; .\;\;\;\;\;\;\;
        \label{g}
\end{eqnarray}
We mention for clarity that the first summand in (\ref{g}) comes 
about as the curvature of the spin connection, whereas the second 
summand arises as the Riemannian curvature of $M$; this can be seen 
more explicitly by writing out ${\nabla}$ with the spin 
derivative $D$ and the Christoffel symbols of the Levi-Civita 
connection on $M$. The third summand in (\ref{g}) is obtained by 
combining (\ref{7}) with (\ref{e}).

With the transformations (\ref{hn}), (\ref{f}), and (\ref{g}), we can 
reduce the third order derivative of the wave function
${\nabla}_j^* {\nabla}^j {\nabla}^k \Psi$ to expressions 
which contain derivatives of $\Psi$ of at most first order.
More precisely, using the Gauss equation,
this allows us to estimate the scalar product
$({\nabla}_j^* {\nabla}^j {\nabla}^k \Psi \:|\: 
{\nabla}_l \Psi)$ in the form
\begin{eqnarray}
({\nabla}_j^* {\nabla}^j {\nabla}^k \Psi \:|\: {\nabla}_l \Psi)
&\leq& C_1 \left( |R_{ijkl}| \:+\: |h_{jk}|^2 \:+\: |{\nabla}_j 
h_{kl}| \right) \:(D\Psi \:|\: D\Psi) \label{h} \\
&&+\:C_2 \left( |{\nabla}_i \bar{R}_{jklm}| \:+\: |h_{ij}| \:|\bar{R}_{klmn}| 
\right) \:(\Psi \:|\: D\Psi)
        \label{i}
\end{eqnarray}
with suitable constants $C_1$ and $C_2$ which are independent of the 
geometry of $M$ and $N$. We multiply both sides of this 
inequality by $\eta$ and integrate over $M$. In the integral over 
(\ref{h}), we estimate with the sup-norm, whereas the integral over 
(\ref{i}) can be bounded using the Schwarz inequality,
\begin{eqnarray*}
\int_M \eta\:({\nabla}_j^* {\nabla}^j {\nabla}^k \Psi \:|\: 
{\nabla}_l \Psi) \:d\mu
&\leq& C_1 \sup_M (\eta \:(|R_{ijkl}| \:+\: |h_{jk}|^2 \:+\: |{\nabla}_j 
h_{kl}|)) \:\int_M (D\Psi \:|\: D\Psi) \:d\mu \\
&&+C_2 \:\sup_M |\Psi| \:
\left\| \eta\: (|{\nabla}_i \bar{R}_{jklm}| + |h_{ij}| \:|\bar{R}_{klmn}|) 
\right\|_{L^2} \| |D\Psi| \|_{L^2} .
\end{eqnarray*}
Finally, we substitute the Positive Energy Theorem (\ref{14b}).
\QED
We now combine the results of Lemma~\ref{lemma1} and~\ref{lemma2}.
\begin{Corollary}
\label{cor1}
There are constants $c_1$ and $c_2$ independent of the geometry of $M$ 
and $N$ such that for any smooth, bounded function $\eta$ with bounded gradient 
on $M$ and the 
Dirac wave function of the Positive Mass Theorem (\ref{14b}),
\begin{eqnarray*}
\lefteqn{ \int_M \eta \left( \sqrt{\bar{R}_{ijkl} \bar{R}^{ijkl}}
- \sqrt{2 \:\bar{R}_{ijk \alpha} \nu^\alpha \:\bar{R}^{ijk \beta} \nu_\beta}
\right)^2 \; (\Psi \:|\: \Psi) \:d\mu } \\
& \leq & c_1 \:(E-|P|) \:\sup_M \left( |\partial_j \eta|
\:|h_{kl}| \:+\: |\eta|\: (|{\nabla}_j h_{kl}| + |R_{ijkl}| +
|h_{ij}|^2) \:+\: |\Delta \eta| \right) \\
&&+ c_2 \:\sqrt{E-|P|} \:\left\|\eta \:(|{\nabla}_k \bar{R}_{ij \alpha
\beta}| \:+\: |h_{ij} \:{\nabla}_k h_{lm}| \:+\: |h_{ij}
\:\bar{R}_{klmn}|) \right\|_{L^2} \; \sup_M |\Psi| \;\;\;.
\end{eqnarray*}
\end{Corollary}
It seems likely to the authors that this inequality is not optimal
in the case of non-zero fundamental form, in the sense that $E-|P|=0$
does not imply that $M^3$ is a submanifold of Minkowski space.
Improvements of the estimate are still under investigation. 
However, in the case $h_{ij} \equiv 0$, the above
inequality is very useful, as we shall see in what follows.

For the rest of this paper we will assume that $h_{ij} \equiv 0$.  Hence, all
of the remaining theorems will concern Riemannian 3-manifolds $(M^3,g)$ which,
by the dominant energy condition (II) and the Gauss equation, must
have non-negative scalar curvature.  Then 
in this zero second fundamental form setting, it follows from the Gauss and 
Codazzi equations
that $\bar{R}_{ijkl} = R_{ijkl}$ and $\bar{R}_{ijk \alpha} \nu^\alpha = 0$,
where $R$ is the Riemannian curvature tensor of $(M^3,g)$.  It also follows
in this setting that the total
momentum is zero, so that the total energy $E$ equals the total mass $m$.

\begin{Corollary}
\label{cor_h=0} There exist positive constants $c_1$ and $c_2$ such that for
any smooth, complete, asymptotically flat manifold $(M^3,g)$ with non-negative
scalar curvature and total mass $m$ and any smooth, bounded function $\eta$
with bounded gradient on $M^3$, 
\begin{eqnarray}
\lefteqn{ \int_M \eta \:R_{ijkl} R^{ijkl} \; (\Psi \:|\: \Psi) \:d\mu }
\nonumber \\
&\leq& c_1 \: m \:\sup_M \left( |\eta| |R_{ijkl}| \:+\: |\Delta \eta| \right)
\:+\: c_2 \:\sqrt{m} \:\left\|\eta \:|{\nabla}_k R_{ij \alpha
\beta}|  \right\|_{L^2} \; \sup_M |\Psi| \;\;\;,
\end{eqnarray}
where $\Psi$ is the 
Dirac wave function of the Positive Mass Theorem (\ref{14b}).
\end{Corollary}
The interesting point of this estimate is that the terms on the right 
side of our inequality all contain factors $m$ or $\sqrt{m}$, 
which vanish when the total mass goes to zero. The disadvantage of 
our estimate is that it involves the Dirac wave function. In order to 
get a more explicit estimate, we shall in the next section derive
a-priori bounds for $\Psi$.

\section{Upper and Lower Bounds for the Norm of the Spinor}
First, we observe that we can use the maximum principle to prove that 
$|\Psi(x)| \le 1$.  To do this, we must derive a second order scalar inequality
for $|\Psi(x)|^2$.
Recall that we are still assuming that $M^3$ has zero second fundamental form,
as we will do for the remainder of the paper.  Then we have
\begin{eqnarray}
\partial_j (\Psi \:|\: \Psi) &=& (D_j \Psi \:|\: \Psi) \:+\: (\Psi
\:|\: D_j \Psi) \label{s1} \\
\Delta (\Psi \:|\: \Psi) &=& {\nabla}_j \left(
(D^j \Psi \:|\: \Psi) \:+\: (\Psi
\:|\: D^j \Psi)  \right) 
\nonumber \\
&=& ({\nabla}_j D^j \Psi \:|\: \Psi) \:+\: 2 \:(D_j \Psi \:|\: 
D^j \Psi) \:+\: (\Psi \:|\: {\nabla}_j D^j \Psi) \nonumber \\
&\stackrel{(\ref{aa}),(\ref{e})}{=}&
-(D^*_j D^j \Psi \:|\: \Psi) \:+\: 2 \:|D \Psi|^2
\:-\: (\Psi \:|\: D^*_j D^j \Psi) \;\;\; . \nonumber 
\end{eqnarray}
Substituting in the Weitzenb\"ock formula (\ref{14b}), we find that, by
(\ref{10}),
\[ \Delta |\Psi|^2 \;\geq\; \frac12 \:(\Psi \:|\: R \:\Psi) \:+\:
2 \:|D \Psi|^2 \;\;\;, \]
where $R$ is the scalar curvature of $(M^3,g)$. Since
$R \geq 0$ in the zero second fundamental form context,
it follows that $|\Psi(x)|^2$ is subharmonic. Using that $|\Psi(x)|^2$ 
goes to one at infinity
by construction, the maximum principle yields that 
\begin{equation}\label{eqn37}
   |\Psi(x)| \le 1
\end{equation}
for all $x$.

To get a lower bound for $|\Psi(x)|$, we let $f(x) = |\Psi(x)|^2$
and observe that
equations (\ref{s1}), (\ref{eqn37}), and (\ref{14b}) imply that 
\[ \int_M |\nabla f|^2 \;\leq\;
 4 \int_M |D \Psi|^2 |\Psi|^2 \;\leq\; 16\pi \, m\;\;\;. \]
Then the Sobolev inequality applied to $1-f(x)$ yields for some positive
constant $c_3$
\begin{equation}
k^2 \left( \int_M (1-f(x))^6 \right)^{1/3}
\;\leq\; c_3 \int_M |\nabla f|^2 \;\leq\; 16\pi \:c_3\: m \;\;\;,
\label{sobolev}
\end{equation}
where $k$ is the isoperimetric constant of $M$ defined to be 
\[ k \;=\; \inf \frac{A}{V^{2/3}} \;\;\;, \]
and the infimum is taken over all smooth regions with volume $V$ and
boundary area $A$ (the Sobolev constant of $M^3$ can be
bounded by $\sqrt{c_3}/k$, where
$c_3$ is a constant independent of the geometry of $M$).
The inequality (\ref{sobolev}) immediately implies
the following lemma, which gives a lower bound for $|\Psi|^2$ except
on a set of small measure.
\begin{Lemma}\label{lem:lower}
Let $k$ be the isoperimetric constant of $M$.  Then for any $c<1$,
\[ |\Psi|^2 \ge c \]
except on a set $D(c)$ with
\begin{equation}
\mbox{Vol}(D(c))^{1/3} \;\leq\;
   \frac{16\pi\:c_3}{(1-c)^2} \: \frac{m}{k^2} \;\;\;.
\end{equation}
\end{Lemma}

\section{Proof of the Main Theorem and Applications}
Our main Theorem~\ref{thm6.2} immediately follows by 
combining the bound (\ref{eqn37}) and Lemma \ref{lem:lower} (where we
set $c=1/2$) with Corollary \ref{cor_h=0}.
We note that the constants $c_1$ ,$c_2$, and $c_3$ could be computed in a
straightforward manner, although we do not carry this out 
since their
actual value is not important for our applications of the theorem.  Also, we 
see that by choosing $\eta$ to be zero everywhere except in a neighborhood of
a given point, Theorem~\ref{thm6.2} yields the following corollary.
\begin{Corollary}
Suppose $\{g_i\}$ is a sequence of smooth, complete, asymptotically flat
metrics on $M^3$ with non-negative scalar curvature and total masses
$\{m_i\}$ which converge to a 
possibly non-smooth limit metric $g$ in the $C^0$ sense.  Let $U$ be the
interior of the set of points where this convergence of metrics is locally 
$C^3$ and nondegenerate.

Then if the metrics $\{g_i\}$ have uniformly positive isoperimetric constants
and their masses $\{m_i\}$ converge to zero, then $g$ is flat in $U$.
\end{Corollary}

Equivalently, we can restate the above corollary in a manner which extends the 
case of equality of the Positive Mass Theorem to manifolds which are not necessarily
smooth.
\begin{Def}
Given a metric $g$ on a manifold $M^3$ which is not necessarily smooth, we say 
that it has {\bf{generalized non-negative scalar curvature}}
if it is the limit in the $C^0$ sense of a sequence
of smooth, complete, asymptotically flat metrics $\{g_i\}$ on 
$M^3$ which have non-negative scalar curvature.  We will also require that $g$
is smooth outside a bounded set, and that its total mass equals the limit of the 
total masses of the smooth metrics.  Furthermore, given such a manifold, let
${\bf{U(M^3)}}$
denote the interior of the set of points in $M^3$ 
where the convergence of metrics is locally $C^3$ and nondegenerate.
\end{Def}

\begin{Thm}\label{last_theorem}
Suppose that $(M^3,g)$ is not necessarily smooth but is complete, asymptotically 
flat, and 
has generalized non-negative scalar curvature, total mass $m$, and isoperimetric 
constant $k$.  Then $m = 0$ and $k > 0$ implies that
$g$ is flat in $U(M^3)$.
\end{Thm}
We note that the above theorem is not true without the requirement that the 
isoperimetric
constant $k>0$.  For example, the induced metric 
on a 3-plane which is tangent to a round 3-sphere in 
Euclidean 4-space has generalized non-negative scalar curvature and is the
limit of portions of space-like Schwarzschild metrics of small mass joined
to the 3-sphere minus a small 
neighborhood of the north pole.  This singular manifold has zero mass, but it is 
not flat everywhere in $U(M^3)$ since $U(M^3)$ equals the whole manifold minus the 
point
where the 3-plane and 3-sphere are tangent.  However, Theorem \ref{last_theorem}
is not contradicted by this example since this singular 
manifold has zero isoperimetric constant.

Among other possible applications,
Theorem \ref{last_theorem} is used in \cite{Bray-Penrose} to
handle the cases of equality of Theorems 1, 9, and 10
of that paper.  These three theorems are
generalizations of the Positive Mass Theorem and give lower bounds on the 
total mass of an asymptotically flat manifold in terms of the geometry of 
the manifold.  In particular,
Theorem 1, which is the main theorem of \cite{Bray-Penrose},
is a slight generalization of the Riemannian Penrose Inequality, which 
states that the total mass of a 3-manifold with non-negative scalar curvature
is greater than or equal to the mass contributed by any black holes it may
contain.  The above Theorem \ref{last_theorem} is then needed to prove that,
if the total mass of the 3-manifold exactly equals the mass contributed by the 
black holes it contains, then the 3-manifold is a Schwarzschild 3-manifold
(defined to be a time-symmetric, space-like slice of the usual 3+1 dimensional
Schwarzschild metric) which corresponds to a single non-rotating black hole
in vacuum.\\[1.5em]
{\em{Acknowledgments:}} We would like to thank Professor S.-T.\ Yau for
many helpful discussions. Also, we want to thank the Harvard Mathematics
Department and the Max Planck Institute for Mathematics in the Sciences,
Leipzig, for their hospitality.

\begin{tabular}{lcl}
\\
Hubert Bray & $\;\;\;\;\;\;\;\;\;$& Felix Finster \\
Department of Mathematics && Max Planck Institute for \\
Massachusetts Institute of Technology && Mathematics in the Sciences \\
77 Massachusetts Avenue && Inselstr.\ 22-26 \\
Cambridge, MA 02139, USA && 04103 Leipzig, Germany \\
\tt{bray@math.mit.edu} && \tt{Felix.Finster@mis.mpg.de}
\end{tabular}

\end{document}